\documentclass[11pt]{amsart}

\usepackage{amsmath,amssymb,amsthm,amsfonts,amscd,flafter,epsf, epsfig,graphicx,verbatim,pinlabel,mathrsfs}
\usepackage[all]{xy}
\usepackage{epsf} 

\title{Comultiplication in link Floer homology and transversely non-simple links}
\author[John A. Baldwin]{John A. Baldwin}
\address {Department of Mathematics, Princeton University\\ Princeton, NJ 08544-1000}
\email {baldwinj@math.princeton.edu}
\thanks{The author was supported by an NSF Postdoctoral Fellowship.}
\date{}

\newcommand\bS{\mathbf{S}}
\newcommand\bx{\mathbf{x}}
\newcommand\by{\mathbf{y}}
\newcommand\bw{\mathbf{w}}

\newcommand\bz{\mathbf{z}}

\newcommand\cfl{\widehat{CFL}}
\newcommand\cflt{\widetilde{CFL}}
\newcommand\cflm{CFL^-}

\newcommand\hfl{\widehat{HFL}}
\newcommand\hflt{\widetilde{HFL}}
\newcommand\hflm{HFL^-}

\newcommand\zzt{\mathbb{Z}_2}

\newtheorem{theorem}{Theorem}[section]
\newtheorem{lemma}[theorem]{Lemma}

\newtheorem{corollary}[theorem]{Corollary}
\newtheorem{proposition}[theorem]{Proposition}

\theoremstyle{definition}

\newtheorem{remark}[theorem]{Remark}

\begin{document}
\begin{abstract} 
For a word $w$ in the braid group $B_n$, we denote by $T_w$ the corresponding transverse braid in $(\mathbb{R}^3,\xi_{rot})$. We exhibit, for any two $g,h \in B_n$, a ``comultiplication'' map on link Floer homology $\widetilde\Phi:\hflt(m(T_{hg}) )\rightarrow \hflt(m(T_g\#T_h))$ which sends $\widetilde\theta(T_{hg})$ to $\widetilde\theta(T_g\#T_h)$. We use this comultiplication map to generate infinitely many new examples of prime topological link types which are not transversely simple.
\end{abstract}

\maketitle

\section{Introduction}

Transverse links feature prominently in the study of contact 3-manifolds. They arise very naturally -- for instance, as binding components of open book decompositions -- and can be used to discern important properties of the contact structures in which they sit (see \cite[Theorem 1.15]{bev} for a recent example). Yet, transverse links, even in the standard tight contact structure, $\xi_{std}$, on $\mathbb{R}^3$, are notoriously difficult to classify up to transverse isotopy. 

A transverse link $T$ comes equipped with two ``classical" invariants which are preserved under transverse isotopy: its topological link type and its self-linking number $sl(T)$. For transverse links with more than one component, it makes sense to refine the notion of self-linking number as follows. Let $T$ and $T'$ be two transverse representatives of some $l$-component topological link type, and suppose there are labelings $T=T_1\cup\dots \cup T_l$ and $T'=T'_1\cup\dots\cup T'_l$ of the components of $T$ and $T'$ such that 
\begin{enumerate}
\item there is a topological isotopy sending $T$ to $T'$ which sends $T_i$ to $T_i'$ for each $i$, and\\
\item $sl(S) = sl(S')$ for any sublinks $S = T_{n_1}\cup\dots \cup T_{n_j}$ and $S' = T_{n_1}'\cup\dots\cup T_{n_j}'.$
\end{enumerate}
Then we say that $T$ and $T'$ have the same \emph{self-linking data}, and we write $\mathcal{SL}(T) = \mathcal{SL}(T')$. A basic question in contact geometry is how to tell, given two transverse representatives, $T$ and $T'$, of some topological link with the same self-linking data, whether $T$ and $T'$ are transversely isotopic; that is, whether the classical data completely determines the \emph{transverse} link type. We say that a topological link type is \emph{transversely simple} if any two transverse representatives $T$ and $T'$ which satisfy $\mathcal{SL}(T) = \mathcal{SL}(T')$ are transversely isotopic. Otherwise, the link type is said to be \emph{transversely non-simple}.

From this point on, we shall restrict our attention to transverse links in the tight rotationally symmetric contact structure, $\xi_{rot},$ on $\mathbb{R}^3,$ which is contactomorphic to $\xi_{std}$. There are several well-known examples of knot types which are transversely simple. Among these are the unknot \cite{yasha5}, torus knots \cite{et4} and the figure eight \cite{EH2}. 


Only recently, however, have knot types been discovered which are not transversely simple. These include a family of 3-braids found by Birman and Menasco \cite{bm3} using the theory of braid foliations; and the (2,3) cable of the (2,3) torus knot, which was shown to be transversely non-simple by Etnyre and Honda using contact-geometric techniques \cite{EH4}. Matsuda and Menasco have since identified two explicit transverse representatives of this cabled torus knot which have identical self-linking numbers, but which are not transversely isotopic \cite{mm}. Their examples take center stage in Section \ref{sec:nonsimple} of this paper.

There has been a flurry of progress in finding transversely non-simple link types in the last couple years, spurred by the discovery of a transverse invariant $\theta$ in link Floer homology by Ozsv{\'a}th, Szab{\'o} and Thurston \cite{oszt}; this discovery, in turn, was made possible by the combinatorial description of $\hflm$ found by Manolescu, Ozsv{\'a}th and Sarkar in \cite{mos} (see also \cite{most}).\footnote{There are several versions of this $\theta$ invariant, denoted by $\theta^-$, $\widehat\theta$ and $\widetilde\theta$. } This $\theta$ invariant is applied by Ng, Ozsv{\'a}th and Thurston in \cite{not} to identify several examples of transversely non-simple links, including the knot $10_{132}$. In \cite{vera}, V{\'e}rtesi proves a connected sum formula for $\theta$, which she weilds to find infinitely many examples of non-prime knots which are transversely non-simple (Kawamura has since proven a similar result without using Floer homology \cite{keiko}; both hers and V{\'e}rtesi's results follow from Etnyre and Honda's work on Legendrian connected sums \cite{EH5}). 

Finding infinite families of transversely non-simple \emph{prime} knots is generally more difficult. Using a slightly different invariant, which we shall denote by $\underline\theta$, derived from knot Floer homology and discovered by Lisca, Ozsv{\'a}th, Szab{\'o} and Stipsicz in \cite{lossz}, Ozsv{\'a}th and Stipsicz identify such an infinite family among two-bridge knots \cite{ost}. And, most recently, Khandhawit and Ng use the invariant $\theta$ to construct a 2-parameter infinite family of prime transversely non-simple knots, which generalizes the example of $10_{132}$ \cite{kng}.

In this paper, we formulate and apply a strategy for generating a slew of new infinite families of transversely non-simple prime links. This strategy hinges on the ``naturality" results below. For a word $w$ in the braid group $B_n$, we denote by $T_w$ the corresponding transverse braid in $(\mathbb{R}^3,\xi_{rot})$. 

\begin{theorem}
\label{thm:nat}
There exists a map on link Floer homology, $$\widetilde\Phi:\hflt(m(T_{w\sigma_i}) )\rightarrow \hflt(m(T_w)),$$ which sends $\widetilde\theta(T_{w\sigma_i})$ to $\widetilde\theta(T_w)$, where $\sigma_i$ is one of the standard generators of $B_n$.
\end{theorem}

This theorem implies the existence of a ``comultiplication" map on link Floer homology, similar in spirit to the map discovered in \cite{bald3}:

\begin{theorem}
\label{thm:comult}
For any two braid words $h$ and $g$ in $B_n$, there exists a map, $$\widetilde\mu:\hflt(m(T_{hg}) )\rightarrow \hflt(m(T_g\#T_h)),$$ which sends $\widetilde\theta(T_{hg})$ to $\widetilde\theta(T_g\#T_h).$
\end{theorem}

One may combine Theorem \ref{thm:comult} with V{\'e}rtesi's result governing the behavior of $\theta$ under connected sums to conclude the following.

\begin{theorem}
\label{thm:nonzero}
If $\widehat\theta(T_g)$ and $\widehat\theta(T_h)$ are both non-zero, then so is $\widehat\theta(T_{hg})$.
\end{theorem}

Here, we sketch one potential way to use these results to find transversely non-simple links. Start with some $w_1$, $w_2 \in B_n$ for which $T_{w_1}$ and $T_{w_2}$ are topologically isotopic and have the same self-linking data, but for which $\widehat\theta(T_{w_1})=0$ while $\widehat\theta(T_{w_2})\neq 0$, so that $T_{w_1}$ and $T_{w_2}$ are not transversely isotopic. Now, choose an $h\in B_n$ for which $\widehat\theta(T_h) \neq 0$. Theorem \ref{thm:nonzero} then implies that $\widehat\theta(T_{hw_2})\neq 0$ as well. If one can show that $\widehat\theta(T_{hw_1}) = 0$, that $T_{hw_1}$ and $T_{hw_2}$ still represent the same topological link type, and that $\mathcal{SL}(T_{hw_1}) = \mathcal{SL}(T_{hw_2})$ (this is automatic if $T_{hw_1}$ and $T_{hw_2}$ are knots), then one may conclude that $T_{hw_1}$ represents a transversely non-simple link type. 

An advantage of this approach for generating new transversely non-simple link types from old over, say, that of \cite{vera, keiko}, is that there is no \emph{a priori} reason to expect that the links so formed are composite. We demonstrate the effectiveness of this approach in Section \ref{sec:nonsimple} of this paper. In doing so, we describe an infinite family of prime transversely non-simple link types (half are knots; the other half are 3-component links) which generalizes the (2,3) cable of the (2,3) torus knot. Moreover, it is clear that this example only scratches the surface of the potential of our more general technique.

Lastly, it is tempting to conjecture that the two invariants $\theta$ and $\underline\theta$ agree for transverse links in $(\mathbb{R}^3,\xi_{rot})$, as they share many formal properties. We prove a partial result in this direction, which follows from Theorem \ref{thm:nat} together with work of Vela-Vick on the $\underline\theta$ invariant \cite{vv}.

\begin{theorem}
\label{thm:equivalent}
$\widehat\theta(T)$ and $\widehat{\underline\theta}(T)$ agree for positive, transverse, connected braids $T$ in $(\mathbb{R}^3,\xi_{rot})$.
\end{theorem}

\subsection*{Organization} In the next section, we outline the relationship between grid diagrams, Legendrian links and their transverse pushoffs. In Section \ref{sec:hfl}, we review the grid diagram construction of link Floer homology and describe some important properties of the transverse invariant $\theta$. In Section \ref{sec:comult}, we prove Theorems \ref{thm:nat}, \ref{thm:comult}, \ref{thm:nonzero} and \ref{thm:equivalent}. And, in Section \ref{sec:nonsimple}, we outline a general strategy for using our comultiplication result to produce new examples of transversely non-simple link types, and we give an infinite family of such examples which are prime. 

\subsection*{Acknowledgements} I wish to thank Lenny Ng for helpful correspondence. His suggestions were key in developing some of the strategy formulated in Section \ref{sec:nonsimple}. Thanks also to the referee for helpful comments.

\section{Grid diagrams, Legendrian and transverse links}

In this section, we provide a brief review of the relationship between Legendrian links in $(\mathbb{R}^3,\xi_{std})$, transverse braids in $(\mathbb{R}^3,\xi_{rot})$ and grid diagrams, largely following the discussion in \cite{kng}. For a more detailed account, see \cite{kng, ngt}. The standard tight contact structure $\xi_{std}$ on $\mathbb{R}^3$ is given as $$\xi_{std} = \text{ker}(dz-ydx).$$ An oriented link $L \subset (\mathbb{R}^3,\xi_{std})$ is called \emph{Legendrian} if it is everywhere tangent to $\xi_{std}$, and \emph{transverse} if it is everywhere transverse to $\xi_{std}$ such that $dz-ydx>0$ along the orientation of $L$. Any smooth link can be perturbed by a $C^0$ isotopy to be Legendrian or transverse. We say that two Legendrian (resp. transverse) links are \emph{Legendrian} (resp. \emph{transversely}) isotopic if they are isotopic through Legendrian (resp. transverse) links.

A Legendrian link $L$ can be perturbed to a transverse link (which is arbitrarily close to $L$ in the $C^{\infty}$ topology) by pushing $L$ along its length in a generic direction transverse to the contact planes in such a way that the orientation of the pushoff agrees with that of $L$. The resulting link $L^+$ is called a \emph{positive transverse pushoff} of $L$. Legendrian isotopic links give rise to transversely isotopic pushoffs. Conversely, every transverse link is the positive transverse pushoff of some Legendrian link; however, two such Legendrian links need not be Legendrian isotopic. The precise relationship between Legendrian and transverse links is best explained via \emph{front projections}.

The front projection of a Legendrian link is its projection onto the $xz$ plane. The front projection of a generic Legendrian link has no vertical tangencies and has only semicubical cusps and transverse double points as singularities. Moreover, at each double point, the slope of the overcrossing is more negative than the slope of the undercrossing. See Figure \ref{fig:example}.c for the front projection of a right-handed Legendrian trefoil. 

\begin{figure}[!htbp]
\labellist 
\hair 2pt 
\tiny
\pinlabel $(a)$ at 135 65
\pinlabel $(b)$ at 273 65
\pinlabel $C$ at 55 15

\endlabellist 
\begin{center}
\includegraphics[height = 1.5cm]{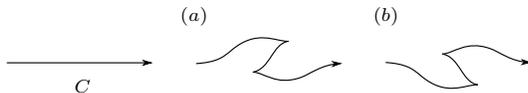}
\caption{\quad \ref{fig:stab}.a and \ref{fig:stab}.b are local pictures of the positive and negative stabilizations, respectively, of a Legendrian link along one of its components $C$.}
\label{fig:stab}
\end{center}
\end{figure}

\begin{figure}[!htbp]
\labellist 
\hair 2pt 
\tiny
\pinlabel $(a)$ at 10 380
\pinlabel $(b)$ at 273 380
\pinlabel $(c)$ at 487 380
\pinlabel $(d)$ at 273 190
\pinlabel $(e)$ at 490 190

\endlabellist 
\begin{center}
\includegraphics[height = 5cm]{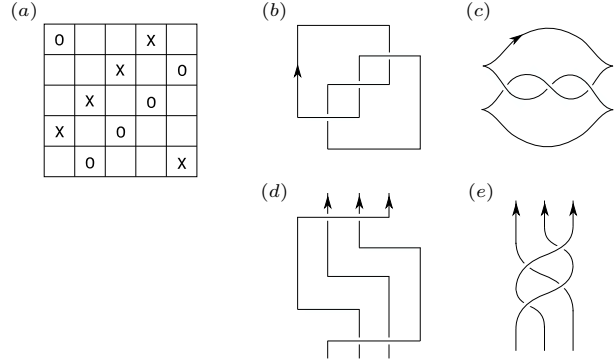}
\caption{\quad In \ref{fig:example}.a, a grid diagram $G$. In \ref{fig:example}.b, the oriented link corresponding to $G$. In \ref{fig:example}.c, a front projection for the Legendrian link $L(G)$. In \ref{fig:example}.d and \ref{fig:example}.e, the braid corresponding to $G$. Here, $w(G) = \sigma_1\sigma_2\sigma_1\sigma_2$.}
\label{fig:example}
\end{center}
\end{figure}

The \emph{positive} (resp. \emph{negative}) \emph{stabilization} of a Legendrian link $L$ along some component $C$ of $L$ is the Legendrian link whose front projection is obtained from that of $L$ by adding a zigzag along $C$ with downward (resp. upward) pointing cusps. See Figure \ref{fig:stab}. Two Legendrian links are said to be \emph{negatively stably isotopic} if they are Legendrian isotopic after each has been negatively stabilized some number of times along some of its components. The following theorem implies that the classification of transverse links up to transverse isotopy is equivalent to the classification of Legendrian links up to Legendrian isotopy and negative stabilization.

\begin{theorem}[\cite{efm,osh}]
Two Legendrian links are negatively stably isotopic if and only if their positive transverse pushoffs are transversely isotopic.
\end{theorem}

Consider the rotationally symmetric tight contact structure on $\mathbb{R}^3$ defined by $$\xi_{rot} = \text{ker}(dz-ydx+xdy).$$ The diffeomorphism of $\mathbb{R}^3$ given by \begin{equation}
\label{eqn:contact}\phi(x,y,z) = (x, 2y, xy+z)\end{equation} sends $\xi_{rot}$ to $\xi_{std}.$ One can define transverse links for $\xi_{rot}$ in the same way that one does for $\xi_{std}$. Since $\phi$ sends a transverse link in $(\mathbb{R}^3,\xi_{rot})$ to a transverse link in $(\mathbb{R}^3,\xi_{std})$, the study of transverse links in $(\mathbb{R}^3,\xi_{std})$ is equivalent to that in $(\mathbb{R}^3,\xi_{rot})$; however, the latter is often more convenient, per the following theorem of Bennequin.

\begin{theorem}[\cite{benn}]
\label{thm:braid2}
Any transverse link in $(\mathbb{R}^3,\xi_{rot})$ is transversely isotopic to a closed braid around the $z$-axis.
\end{theorem}
 
Theorem \ref{thm:braid2} allows us to use braid-theoretic techniques to study transverse links. For a braid word $w\in B_n$, we let $T_w$ denote the corresponding transverse braid around the $z$-axis. Braid words which are conjugate in $B_n$ clearly correspond to transversely isotopic links. Recall that, for $w\in B_n$, a \emph{positive} (resp. \emph{negative}) \emph{braid stabilization} of $w$ is the operation which replaces $w$ by the word $w\sigma_n$ (resp. $w\sigma_n^{-1}$) in $ B_{n+1}$. We will also refer to $T_{w\sigma_n}$ (resp. $T_{w\sigma_n^{-1}}$) as the \emph{positive} (resp. \emph{negative}) \emph{braid stabilization} of the transverse link $T_{w}$. The following theorem makes precise the relationship between braids and transverse links in $(\mathbb{R}^3,\xi_{rot})$.

\begin{theorem}[\cite{osh,wrinkle}]
\label{thm:markov}
For $w\in B_n$ and $w' \in B_m$, the transverse links $T_{w}$ and $T_{w'}$ are transversely isotopic in $(\mathbb{R}^3,\xi_{rot})$ if and only if $w$ and $w'$ are related by a sequence of conjugations and positive braid stabilizations and destabilizations.
\end{theorem}

In Section \ref{sec:nonsimple}, we use a braid operation called an \emph{exchange move}. If $a$, $b$ and $c$ in $B_n$ are words in the generators $\sigma_2,\dots,\sigma_{n-1}$, then an exchange move is the operation which replaces the word $w_1 = a\sigma_1b\sigma_1^{-1}c$ with the word $w_2 = a\sigma_1^{-1}b\sigma_1c$. An exchange move is actually just a composition of conjugations, one positive braid stabilization and one positive destabilization, and so the link $T_{w_1}$ is transversely isotopic to $T_{w_2}$ (see, for example, \cite{ngt}).

It bears mentioning that the self-linking number of a transverse link admits a nice formulation in the language of braids. If $\Sigma$ is a Seifert surface for a transverse link $T$, then the vector bundle $\xi_{rot}|_{\Sigma}$ is trivial and, therefore, has a non-zero section $v$. Recall that the \emph{self-linking number} of $T$ is defined by $$sl(T)=lk(T,T'),$$ where $T'$ is a pushoff of $T$ in the direction of $v$. Any two links which are transversely isotopic have identical self-linking numbers. For a word $w\in B_n$, the self-linking number of $T_w$ is given simply by $a(w)-n$, where $a(w)$ is the algebraic length of $w$.

In what remains of this section, we describe a relationship between the front diagram of a Legendrian link in $(\mathbb{R}^3,\xi_{std})$ and a braid representation of its positive transverse pushoff, thought of as a transverse link in $(\mathbb{R}^3,\xi_{rot})$. Grid diagrams provide the necessary connection.

A \emph{grid diagram} $G$ is an $k \times k$ square grid along with a collection of $k$ $X$'s and $k$ $O$'s contained in these squares such that every row and column contains exactly one $O$ and one $X$ and no square contains both an $O$ and an $X$. See Figure \ref{fig:example}.a. We call $k$ the \emph{grid number} of $G$. One can produce an oriented link diagram $L$ from $G$ by drawing a horizontal segment from the $O$'s to the $X$'s in each row and a vertical segment from the $X$'s to the $O$'s in each column so that the horizontal segments pass over the vertical segments (this is the convention used in \cite{kng}, and the opposite of the convention in \cite{mos}; see \cite{ngt} for a discussion on the relationship between the two conventions), as in Figure \ref{fig:example}.b. By rotating $L$ $45^\circ$ clockwise, and then smoothing the upward and downward pointing corners and turning the leftward and rightward pointing corners into cusps, one obtains the front projection of a Legendrian link, as in Figure \ref{fig:example}.c. Let us denote this Legendrian link by $L(G)$.

Alternatively, one can construct a braid diagram from $G$ by drawing a horizontal segment from the $O$'s to the $X$'s in each row, as before, and drawing a vertical segment from the $X$'s to the $O$'s for each column in which the marking $X$ lies under the marking $O$. For those columns in which the $X$ is above the $O$, we draw two vertical segments: one from the $X$ up to the top of the grid diagram, and the other from the bottom of the grid diagram up to the $O$. As before, we require that the horizontal segments pass over the vertical segments. Note that all vertical segments are oriented upwards and that the closure of the diagram we have constructed is a braid. See Figures \ref{fig:example}.d and \ref{fig:example}.e for an example of this procedure. Let us denote the corresponding braid word by $w(G)$, read from the bottom up. The relationship between $T_{w(G)}$ and $(L(G))^+$ is expressed in the proposition below.

\begin{proposition}[{\rm \cite[Proposition 3]{kng}}]
The contactomorphism $\phi$ from $(\mathbb{R}^3,\xi_{rot})$ to $(\mathbb{R}^3,\xi_{std})$ defined in Equation (\ref{eqn:contact}) sends the transverse link $T_{w(G)}$ to a link which is transversely isotopic to $(L(G))^+$.
\end{proposition}

\section{Link Floer homology and the transverse invariant}
\label{sec:hfl}

In this section, we describe the grid diagram formulation of link Floer homology discovered in \cite{mos,most}. Let $G$ be a grid diagram for a link $L$ and suppose that $G$ has grid number $k$. From this point forward, we think of $G$ as a \emph{toroidal} grid diagram -- that is, we identify the top and bottom sides of $G$ and the right and left sides of $G$ -- so that the horizontal and vertical lines become $k$ horizontal and $k$ vertical circles. Let $\mathbb{O}$ and $\mathbb{X}$ denote the sets of markings $\{O_i\}_{i=1}^k$ and $\{X_i\}_{i=1}^k$, respectively.

We associate to $G$ a chain complex $(\cflm(m(L)), \partial^-)$ as follows. The generators of $\cflm(L)$ are one-to-one correspondences between the horizontal and vertical circles of $G$. Equivalently, we may think of a generator as a set of $k$ intersection points between the horizontal and vertical circles, such that no intersection point appears on more than one horizontal circle or on more than one vertical circle. We denote this set of generators by $\mathbf{S}(G)$. Then, $\cflm(m(L))$ is defined to be the free $\zzt[U_1,\dots,U_k]$-module generated by the elements of $\mathbf{S}(G)$, where the $U_i$ are formal variables corresponding to the markings $O_i$.

\begin{figure}[!htbp]
\labellist 
\hair 2pt 
\tiny

\endlabellist 
\begin{center}
\includegraphics[height = 3cm]{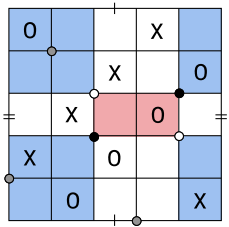}
\caption{\quad A grid diagram $G$ for the right-handed trefoil. The generator $\bx$ comprises the black and gray intersection points while $\by$ comprises the white and gray intersection points. $Rect_G(\bx, \by)$ contains the shaded rectangles in red and blue, while $Rect_G^o(\bx,\by)$ contains only the red rectangle.}
\label{fig:rect}
\end{center}
\end{figure}

For $\bx, \by \in \bS(G)$, we let $Rect_G(\bx,\by)$ denote the space of embedded rectangles in $G$ with the following properties. $Rect_G(\bx,\by)$ is empty unless $\bx$ and $\by$ coincide at $k-2$ points. An element $r\in Rect_G(\bx,\by)$ is an embedded disk on the toroidal grid $G$ whose edges are arcs on the horizontal and vertical circles and whose four corners are intersection points in $\bx \,\cup\, \by$. Moreover, we stipulate that if we traverse each horizontal boundary component of $r$ in the direction specified by the induced orientation on $\partial r$, then this horizontal arc is oriented from a point in $\bx$ to a point in $\by$. If $Rect_G(\bx,\by)$ is non-empty, then it consists of exactly two rectangles. See Figure \ref{fig:rect} for an example. We let $Rect_G^o(\bx,\by)$ denote the space of $r\in Rect_G^o(\bx, \by)$ for which $r\cap\mathbb{X}= \text{Int}(r)\cap \mathbf{x} = \emptyset$.

The module $\cflm(m(L))$ is endowed with an endomorphism $$\partial^-: \cflm(m(L))\rightarrow \cflm(m(L)),$$ defined on $\bS(G)$ by $$\partial^-(\bx) = \sum_{\by\in\bS(G)}\,\,\sum_{r\in Rect_{G}^o(\bx,\by)} U_1^{O_1(r)}\cdots U_k^{O_n(r)} \cdot \by.$$ Here, $O_i(r)$ denotes the number of times the marking $O_i$ appears in $r$. The map $\partial^-$ is a differential, and, so, gives rise to a chain complex $(\cflm(m(L)), \partial^-).$ The homology of this chain complex, $\hflm(m(L)) = H_*(\cflm(m(L)), \partial^-)$, is an invariant of the link $L$, and agrees with the \emph{link Floer homology} of $m(L)$ defined in \cite{osz19}. It bears mentioning that the complex $(\cflm(m(L),\partial^-)$ comes equipped with \emph{Maslov} and \emph{Alexander} gradings, which are then inherited by $\hflm(m(L))$; however, we will not discuss these gradings further as they play no role in this paper.

Suppose that the link $L$ has $l$ components. If $O_i$ and $O_j$ lie on the same component of $L$, then multiplication by $U_i$ in $(\cflm(m(L)), \partial^-)$ is chain homotopic to multiplication by $U_j$, and, so, these multiplications induce the same maps on $\hflm(m(L))$ \cite[Lemma 2.9]{most}. So, if we label the markings in $\mathbb{O}$ so that $O_1,\dots,O_l$ lie on different components, then we can think of $\hflm(m(L))$ as a module over $\zzt[U_1,\dots,U_l]$. 

Setting $U_1 = \dots = U_l=0$, one obtains a chain complex $(\cfl(m(L)), \widehat\partial)$ whose homology we denote by $\hfl(m(L))$. The latter is a bi-graded vector space over $\zzt$, whose graded Euler characteristic is some normalization of the multivariable Alexander polynomial of $m(L)$ \cite{osz19}. If one sets $U_1 = \dots = U_k = 0$, one obtains a chain complex $(\cflt(m(L)), \widetilde\partial)$ whose homology we denote by $\hflt(m(L))$. The group $\hfl(m(L))$ determines $\hflt(m(L))$. Specifically, if we let $n_i$, for $i=1,\dots,l$, denote the number of markings in $\mathbb{O}$ on the ith component of $L$, then $$\hflt(m(L)) = \hfl(m(L)) \otimes \bigotimes_{i=1}^l V_i^{\otimes (n_i-1)},$$ where $V_i$ is a fixed two dimensional vector space \cite[Proposition 2.13]{most}, and the quotient map $$j:\cfl(m(L))\rightarrow \cflt(m(L))$$ induces an injection $j_*$ on homology.

The element $\bz^+(G) \in \bS(G)$, which consists of the intersection points at the upper right corners of the squares in $G$ containing the markings in $\mathbb{X}$, is clearly a cycle in $(\cflm(m(L)), \partial^-)$ (and, hence, also in the other chain complexes). If $T$ is the transverse link in $(\mathbb{R}^3,\xi_{rot})$ corresponding to the braid obtained from $G$ as in Figure \ref{fig:example}.e, then $T$ is topologically isotopic to $L$, and the image of $\bz^+(G)$ in $\hflm(m(T))$ is the transverse invariant $\theta^-(T)$ defined in \cite{oszt}. The images of $\bz^+(G)$ in $\hfl(m(T))$ and $\hflt(m(T))$ are likewise denoted $\widehat\theta(T)$ and $\widetilde\theta(T)$, and are invariants of the transverse link $T$ as well. Moreover, the map $j_*$ sends $\widehat\theta(T)$ to $\widetilde\theta(T)$; in particular, $\widehat\theta(T)=0$ if and only if $\widetilde\theta(T)=0$. The theorem below makes these statements about invariance precise.

\begin{theorem}[{\rm \cite[Theorem 7.1]{oszt}}]
\label{thm:invt}
Suppose that $G$ and $G'$ are two grid diagrams whose associated braids $T$ and $T'$ are transversely isotopic. Then, there is an isomorphism $$f^o_*:HFL^o(m(T))\rightarrow HFL^o(m(T')),$$ induced by a chain map $f^o$, which sends $\theta^o(T)$ to $\theta^o(T')$. 
\end{theorem}

Here, the superscript ``\,$^o$\," is meant to indicate that this theorem holds for any of the three versions of link Floer homology described above. In particular, if $T$ and $T'$ are two transverse links for which $\widehat\theta(T) \neq 0$ and $\widehat\theta(T') = 0$, then $T$ and $T'$ are not transversely isotopic (the invariant $\theta^-(T)$ is always non-zero and non-$U_i$-torsion in $\hflm(m(T))$ \cite[Theorem 7.3]{oszt}). These transverse invariants also behave nicely under negative braid stabilizations.

\begin{theorem}[{\rm \cite[Theorem 7.2]{oszt}}]
\label{thm:stab}
Suppose that $G$ and $G'$ are two grid diagrams with associated braids $T$ and $T'$, and suppose that $T'$ is obtained from $T$ by performing a negative braid stabilization along the ith component of $T$. Then, there is an isomorphism $$f^-_*:\hflm(m(T))\rightarrow \hflm(m(T')),$$ induced by a chain map $f^-$, which sends $\theta(T)$ to $U_i\cdot\theta(T')$.
\end{theorem}

Since multiplication by $U_i$ is the same as multiplication by zero on $\hfl$ and $\hflt$, we obtain the following corollary.

\begin{corollary}
\label{cor:stab}
If $T'$ is obtained from a transverse braid $T$ by performing a negative braid stabilization along some component of $T$, then $\widehat\theta(T') = \widetilde\theta(T')=0.$
\end{corollary}

\section{The map $\Phi$ and comultiplication}
\label{sec:comult}

\begin{figure}[!htbp]
\labellist 
\hair 2pt 
\tiny

\pinlabel $T_{w\sigma_i}$ at 440 50
\pinlabel $T_{w}$ at 600 50
\tiny
\pinlabel $\beta$ at -7 107
\pinlabel $\gamma$ at -7 137

\endlabellist 
\begin{center}
\includegraphics[height = 5cm]{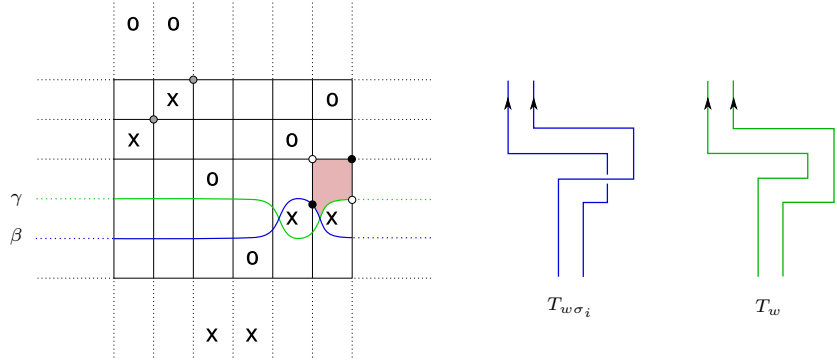}
\caption{\quad A portion of the grid diagrams $G_{\beta}$ and $G_{\gamma}$. The cycle $\mathbf{z}^+(G_{\beta})$ is shown as a collection of black and gray dots, while $\mathbf{z}^+(G_{\gamma})$ is represented by the white and gray dots.}
\label{fig:pentagon}
\end{center}
\end{figure}

Fix some $w\in B_n$ and some $i\in\{1,\dots,n-1\}$. Figure \ref{fig:pentagon} shows simultaneously a portion of a grid diagram $G_{\beta}$ for $T_{w\sigma_i}$ and the corresponding portion of a grid diagram $G_{\gamma}$ for $T_{w}$. The grid diagrams $G_{\beta}$ and $G_{\gamma}$ are the same except that $G_{\beta}$ uses the horizontal curve $\beta$ while $G_{\gamma}$ uses the horizontal curve $\gamma$. Let $k$ denote their common grid number.

For $\bx \in \bS(G_{\beta})$ and $\by \in \bS(G_{\gamma})$, let $Pent_{\beta\gamma}(\bx,\by)$ denote the space of embedded pentagons with the following properties. $Pent_{\beta\gamma}(\bx,\by)$ is empty unless $\bx$ and $\by$ coincide at $k-2$ points. An element $p\in Pent_{\beta\gamma}(\bx,\by)$ is an embedded disk in the torus whose boundary consists of five arcs, each contained in horizontal or vertical circles. We stipulate that under the orientation induced on the boundary of $p$, the boundary may be traversed as follows. Start at the component of $\bx$ on the curve $\beta$ and proceed along an arc contained in $\beta$ until we arrive at the right-most intersection point between $\beta$ and $\gamma$; next, proceed along an arc contained in $\gamma$ until we reach the component of $\by$ contained in $\gamma$; next, follow the arc of a vertical circle until we arrive at a component of $\bx$; then, proceed along the arc of a horizontal circle until we arrive at a component of $\by$; finally, follow an arc contained in a vertical circle back to the initial component of $\bx$. Let $Pent_{\beta\gamma}^o(\bx,\by)$ denote the space of $p\in Pent_{\beta\gamma}(\bx,\by)$ for which $p\cap\mathbb{X}= \text{Int}(p)\cap \mathbf{x} = \emptyset$.

We construct a map $$\phi^-:\cflm(m(T_{w\sigma_i}))\rightarrow \cflm(m(T_w))$$ of $\zzt[U_1,\dots,U_k]$-modules as follows. For $\bx \in \bS(G_{\beta})$, let $$\phi^-(\bx) = \sum_{\by\in\bS(G_{\gamma})}\,\,\sum_{p\in Pent_{\beta\gamma}^o(\bx,\by)} U_1^{O_1(p)}\cdots U_k^{O_n(p)} \cdot \by.$$ We then define $$\widetilde\phi:\cflt(m(T_{w\sigma_i}))\rightarrow \cflt(m(T_w))$$ to be the map on $\cflt$ induced by $\phi^-$. In other words, $\widetilde\phi$ counts pentagons in $Pent_{\beta\gamma}^o(\bx,\by)$ which also miss the $\mathbb{O}$ basepoints. (This construction is inspired by the proof of commutation invariance in \cite[Section 3.1]{most}.) 

\begin{remark}
Unlike $\widetilde\phi$, the map $\phi^-$ is not necessarily a chain map.
\end{remark}

The juxtaposition $p*r$ of any $p\in Pent^o_{\beta\gamma}(\bx,\by)$ with any rectangle $r \in Rect^o_{G_\gamma}(\by,\bw)$ such that $p\cap\mathbb{O}=r\cap\mathbb{O}=\emptyset$ has precisely one such decomposition and exactly one other decomposition as $r'*p'$, where $r' \in Rect^o_{G_\beta}(\bx,\by')$ and $p'\in Pent^o_{\beta\gamma}(\by',\bw)$ and $r'\cap\mathbb{O}=p'\cap\mathbb{O}=\emptyset$. It follows that $\widetilde\phi$ is a chain map and, so, induces a map $$\widetilde\Phi:\hflt(m(T_{w\sigma_i}))\rightarrow \hflt(m(T_w)).$$ Moreover, it is clear that $Pent^o_{\beta\gamma}(\mathbf{z}^+(G_{\beta}), \mathbf{z}^+(G_{\gamma}))$ consists only of the shaded pentagon shown in Figure \ref{fig:pentagon}, and that $Pent^o_{\beta\gamma}(\mathbf{z}^+(G_{\beta}), \mathbf{y})$ is empty for $\mathbf{y} \neq  \mathbf{z}^+(G_{\gamma}).$ Therefore, $\widetilde\Phi$ sends $\widetilde\theta(T_{w\sigma_i})$ to $\widetilde\theta(T_w)$, proving Theorem \ref{thm:nat}.

\begin{figure}[!htbp]
\labellist 
\hair 2pt 
\tiny
\pinlabel $(a)$ at 30 850
\pinlabel $(b)$ at 393 850
\pinlabel $(c)$ at 900 850
\pinlabel $(d)$ at 1400 850
\pinlabel $(e)$ at 1900 850

\small
\pinlabel $h$ at 160 398
\pinlabel $h$ at 515 398
\pinlabel $h$ at 1020 398
\pinlabel $\reflectbox{h}$ at 1520 398
\pinlabel $\reflectbox{h}$ at 2023 397
\pinlabel $g$ at 160 737
\pinlabel $g$ at 515 737
\pinlabel $g$ at 1020 737
\pinlabel $g$ at 1680 737
\pinlabel $g$ at 2190 737
\endlabellist 
\begin{center}
\includegraphics[height = 5cm]{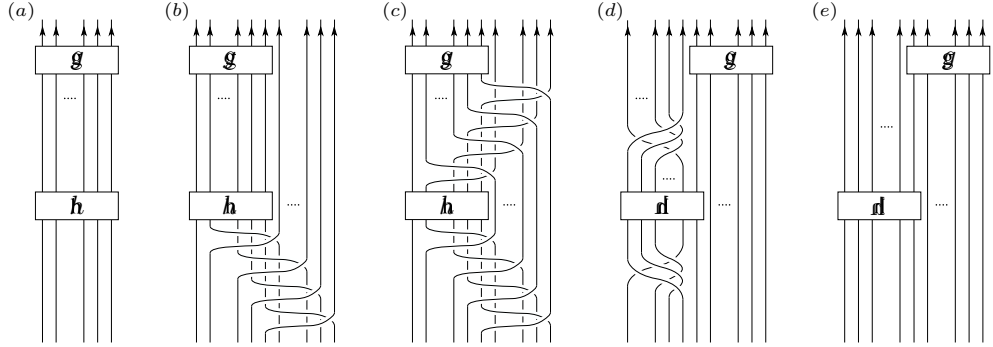}
\caption{\quad $T_{hg}$ is transversely isotopic to the braid in \ref{fig:braid}.b, which, after introducing negative crossings (or, equivalently, up to braid isotopy, getting rid of positive crossings) is transversely isotopic to $T_g \#T_h,$ which is represented by the braid in \ref{fig:braid}.e. }
\label{fig:braid}
\end{center}
\end{figure}

The more general comultiplication fact stated in Theorem \ref{thm:comult} follows from the above result together with the sequence of braid moves depicted in Figure \ref{fig:braid}. The braid in Figure \ref{fig:braid}.a represents the transverse link $T_{hg}$. The braid in \ref{fig:braid}.b is obtained from that in \ref{fig:braid}.a by a mixture of isotopy and positive stabilizations. The braid in \ref{fig:braid}.c is obtained from that in \ref{fig:braid}.b by isotopy followed by the introduction of negative crossings. The braid in \ref{fig:braid}.e is isotopic to the braids in \ref{fig:braid}.c and \ref{fig:braid}.d, and represents the connected sum of the transverse links $T_g$ and $T_h$ (for the latter statement, see \cite{bm2}). Therefore, a composition of the maps $\widetilde\Phi$ described above (one for each negative crossing introduced in going from \ref{fig:braid}.b to \ref{fig:braid}.c) yields a map $$\widetilde\mu:\hflt(m(T_{hg}) )\rightarrow \hflt(m(T_g\#T_h))$$ which sends $\widetilde\theta(T_{hg})$ to $\widetilde\theta(T_g\#T_h).$ 

Suppose that $T_g \# T_h$ is any connected sum of $T_g$ and $T_h$. In \cite{vera}, V{\'e}rtesi proves the following refinement of the Kunneth formula described in \cite[Theorem 1.4]{osz19}. (Her proof is actually for the analogous result in knot Floer homology, but it extends in an obvious manner to a proof of the theorem below.)


\begin{theorem}
\label{thm:connectedsum}
There is an isomorphism, $$\hfl( m(T_g \# T_h)) \cong \hfl( m(T_g))\otimes_{\zzt} \hfl(m(T_h)),$$ under which $\widehat\theta(T_g \# T_h)$ is identified with $\widehat\theta(T_g)\otimes \widehat\theta(T_h)$.
\end{theorem}

V{\'e}rtesi's theorem, used in combination with the comultiplication map $\widetilde{\mu}$, may be applied to prove Theorem \ref{thm:nonzero}.

\begin{proof}[Proof of Theorem \ref{thm:nonzero}]
Recall from the previous section that $\widehat\theta(T_w)$ is non-zero if and only if $\widetilde\theta(T_w)$ is non-zero. If $\widehat\theta(T_g)$ and $\widehat\theta(T_h)$ are both non-zero, then, by Theorem \ref{thm:connectedsum}, so is $\widehat\theta(T_g \# T_h)$, and, hence, so is $\widetilde\theta(T_g \# T_h)$. Since $\widetilde\mu$ sends $\widetilde\theta(T_{hg})$ to $\widetilde\theta(T_g\#T_h),$ this implies that $\widetilde\theta(T_{hg})$ is non-zero, and, hence, so is $\widehat\theta(T_{hg})$.
\end{proof}

Recall that a braid $T_g$ is said to be \emph{quasipositive} if $g \in B_n$ can be expressed as a product of conjugates of the form $w\sigma_iw^{-1}$, where $w$ is any word in $B_n$.

\begin{corollary}
\label{cor:qp}
If $T_g$ is a quasipositive braid, then $\widehat\theta(T_g)\neq0.$
\end{corollary}

\begin{figure}[!htbp]
\labellist 
\hair 2pt 
\small
\pinlabel $n$ at 242 105

\endlabellist 
\begin{center}
\includegraphics[height = 3cm]{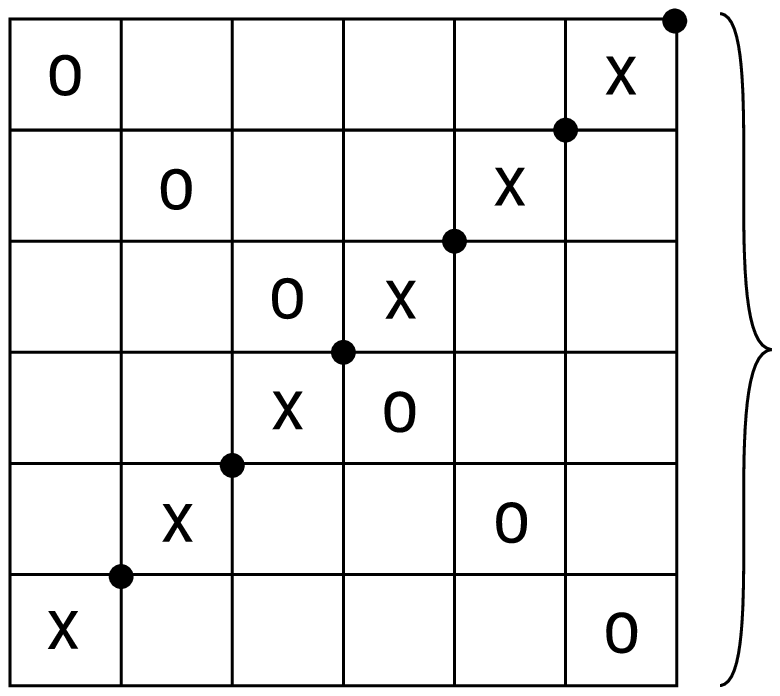}
\caption{\quad A grid diagram $G_n$ for the trivial braid $I_n$. It is straightforward to check that the cycle $\mathbf{z}^+(G_{n}) \in \mathbf{S}(G_n)$, represented by the collection of black dots, is not a boundary in $\cflt(m(I_n))$. }
\label{fig:trivial}
\end{center}
\end{figure}

\begin{proof}[Proof of Corollary \ref{cor:qp}]
If $g$ is a product of $m$ conjugates as above, then after resolving the corresponding $m$ positive crossings, one obtains a braid isotopic to $I_n$, the trivial $n$-braid. Therefore, a composition of $m$ of the maps $\widetilde\Phi$ sends $\widetilde\theta(T_g)$ to $\widetilde\theta(I_n)$. Moreover, one sees by glancing at the grid diagram for $I_n$ in Figure \ref{fig:trivial} that $\widetilde\theta(T_{I_n})\neq 0$. Therefore, $\widetilde\theta(T_g) \neq 0$ and the same is true of $\widehat\theta(T_{g})$.
\end{proof}

\begin{proof}[Proof of Theorem \ref{thm:equivalent}]
Suppose that $T_w$ is a positive braid with one component. Then $\widehat\theta(T_w) \neq 0$, by the corollary above; also, $T_w$ is a fibered knot \cite{bw}. Moreover, $\widehat\theta(T_w)$ lies in Alexander grading $(sl(T_w)+1)/2$, which, in this case, is simply the genus of $T_w$ \cite{oszt}. Therefore, $\widehat\theta(T_w)$ is the unique generator of $\hfl(T_w,g(T_w))$. To show that $\widehat\theta(T_w) = \widehat{\underline\theta}(T_w)$, it suffices to prove that $\widehat{\underline\theta}(T_w)$ is non-zero as well. Fortunately, this has been shown by Vela-Vick in \cite{vv}.
\end{proof}

\section{Finding new transversely non-simple links}
\label{sec:nonsimple}

In this section, we outline and apply one strategy for using comultiplication (in particular, Theorem \ref{thm:nonzero}) to generate a plethora of new examples of transversely non-simple link types. Consider the braid words $$w_1 = a\sigma_1^{m}b\sigma_1^{-1}c\,\,\,\,\,\,\, and \,\,\,\,\,\,\,w_2 = a\sigma_1^{-1}b\sigma_1^{m}c$$ in $B_n$, where $a$, $b$ and $c$ are words in the generators $\sigma_2,\dots,\sigma_{n-1}$. The transverse braids $T_{w_1}$ and $T_{w_2}$ are said to be related by a \emph{negative flype} and, in particular, represent the same topological link type. If, in addition, $m$ is odd, or if $m$ is even and the two strands which cross according to $\sigma_1^m$ belong to the same component of $T_{w_1}$, then $\mathcal{SL}(T_{w_1}) = \mathcal{SL}(T_{w_2})$.

Suppose that $\widehat\theta(T_{w_1})=0$ and $\widehat\theta(T_{w_2}) \neq 0$. The idea is to find a word $h$ in the generators $\sigma_2,\dots,\sigma_{n-1}$ with $\widehat\theta(T_h) \neq 0.$ Theorem \ref{thm:nonzero} would then imply that $\widehat\theta(T_{hw_2}) \neq 0$. If it is also true that $\widehat\theta(T_{hw_1})=0$, then $T_{hw_1}$ and $T_{hw_2}$ are not transversely isotopic although they are topologically isotopic. We would like to find examples which also satisfy $\mathcal{SL}(T_{hw_1}) = \mathcal{SL}(T_{hw_2})$ (if $T_{hw_1}$ is a knot, this is automatic) so as to produce topological link types which are not transversely simple. One nice feature of this proposed method, which differs from that in \cite{vera}, is that there is no reason to believe \emph{a priori} that the link $T_{hw_1}$ so obtained is composite.

In principle, Theorem \ref{thm:nonzero} eliminates half of the work in this scenario -- namely, showing that $\widehat\theta(T_{hw_2}) \neq 0$. In practice, one would like to find examples in which the other half -- showing that $\widehat\theta(T_{hw_1})$ is zero -- is very easy. To that end, one strategy is to pick an example in which $T_{w_1}$ is transversely isotopic to a braid which can be negatively destabilized, and to show that the same is true of the braid $T_{hw_1}$, which would guarantee that $\widehat\theta(T_{hw_1})=0$ by Corollary \ref{cor:stab}.  In particular, $T_{w_1}$ must belong to a topological link type with a transverse representative (that is, $T_{w_2}$) which does not maximize self-linking number, but which cannot be negatively destabilized. 

The most well-known such link type is that of the $(2,3)$ cable of the $(2,3)$ torus knot. In \cite{EH4}, Etnyre and Honda prove the following.

\begin{proposition}
\label{prop:EH}
The $(2,3)$ cable of the $(2,3)$ torus knot has two Legendrian representatives, $L_1$ and $L_2$, both with $tb=5$ and $r=2$, for which $L_1$ is the positive (Legendrian) stabilization of a Legendrian knot while $L_2$ is not. Moreover, $L_1$ and $L_2$ are not Legendrian isotopic after any number of negative (Legendian) stabilizations.
\end{proposition}

That $L_1$ and $L_2$ are not Legendrian isotopic after any number of negative stabilizations implies that their transverse pushoffs, $L_1^+$ and $L_2^+$, are not transversely isotopic (yet, they both have $sl = 3$). Moreover, since $L_1$ is the positive stabilization of a Legendrian knot, its pushoff $L_1^+$ is transversely isotopic to the negative stabilization of some transverse braid. 

\begin{figure}[!htbp]
\labellist 
\hair 2pt 
\tiny
\pinlabel $(a)$ at 60 817
\pinlabel $(b)$ at 390 817
\pinlabel $(c)$ at 732 817
\pinlabel $(d)$ at 1087 817
\pinlabel $(e)$ at 1450 817

\pinlabel $(a')$ at 60 385
\pinlabel $(b')$ at 390 385
\pinlabel $(c')$ at 736 385
\pinlabel $(d')$ at 1100 385
\pinlabel $(e')$ at 1465 385

\endlabellist 
\begin{center}
\includegraphics[height = 6.9cm]{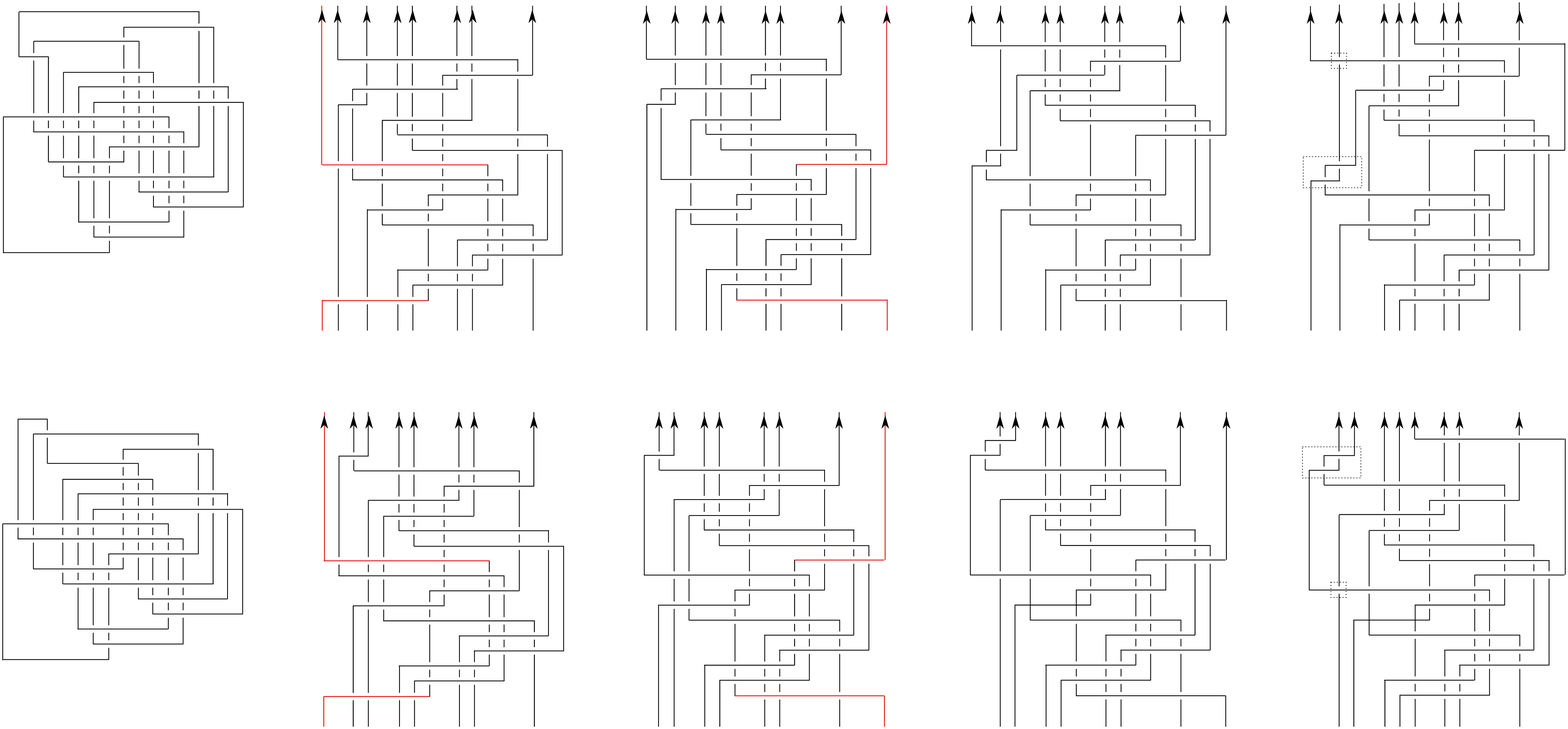}
\caption{\quad On the top, the rectangular diagrams for $L_1$ and its transverse pushoff $L_1^+$. On the bottom, those for $L_2$ and $L_2^+$. The circled regions in \ref{fig:grids}.e and \ref{fig:grids}.e$'$ indicate that $L_1^+$ and $L_2^+$ are transversely isotopic to braids related by a negative flype.}
\label{fig:grids}
\end{center}
\end{figure}

Matsuda and Menasco have since given explicit forms for $L_1$ and $L_2$ \cite{mm}. Figures \ref{fig:grids}.a and \ref{fig:grids}.a$'$ depict the rectangular diagrams corresponding to slightly modified versions of these forms (ours are derived from the front diagrams in \cite[Figure 6]{not}). Figures \ref{fig:grids}.b and \ref{fig:grids}.b$'$ show the rectangular braid diagrams for the transverse pushoffs $L_1^+$ and $L_2^+$, respectively. The braids in \ref{fig:grids}.c and \ref{fig:grids}.c$'$ are obtained from those in \ref{fig:grids}.b and \ref{fig:grids}.b$'$ by isotoping the red arcs as indicated, and the braids in \ref{fig:grids}.d and \ref{fig:grids}.d$'$ are obtained from those in \ref{fig:grids}.c and \ref{fig:grids}.c$'$ after additional simple isotopies and conjugations. The braids in \ref{fig:grids}.e and \ref{fig:grids}.e$'$ are obtained from those in \ref{fig:grids}.d and \ref{fig:grids}.d$'$ by conjugation, and they are related to one another by a negative flype. Indeed, Figure \ref{fig:grids} shows that $L_1^+$ and $L_2^+$ are transversely isotopic to the transverse braids $T_{w_1}$ and $T_{w_2}$, respectively, where $w_1 = a\sigma_1^2 b\sigma_1^{-1}c$, $w_2 = a\sigma_1^{-1} b\sigma_1^2c$, \begin{eqnarray*}a&=&\sigma_4\sigma_3\sigma_5\sigma_6\sigma_4\sigma_5\sigma_5\sigma_6\sigma_4\sigma_5\sigma_7\sigma_6\sigma_5^{-1}\sigma_4^{-1}\sigma_3^{-1}\sigma_2\sigma_3\sigma_3\sigma_4\sigma_5\sigma_4^{-1}\sigma_3^{-1}\sigma_2^{-1},  \\ 
b&=& \sigma_5\sigma_6\sigma_7\sigma_6^{-1}\sigma_5^{-1}\sigma_4^{-1}\sigma_6^{-1}\sigma_5^{-1}\sigma_4^{-1}\sigma_3\sigma_4\sigma_5\sigma_2\sigma_3\sigma_4\sigma_4\sigma_5\sigma_6\sigma_5^{-1}\sigma_4^{-1}\sigma_3^{-1}\sigma_2^{-1}, \text{ and}\\
c&=&\sigma_7^{-1}\sigma_6^{-1}\sigma_5^{-1}.
\end{eqnarray*} 

\begin{figure}[!htbp]
\labellist 
\hair 2pt 
\tiny
\pinlabel $(a)$ at 15 915
\pinlabel $(b)$ at 370 915
\pinlabel $(c)$ at 727 915
\pinlabel $(d)$ at 1060 915
\pinlabel $(e)$ at 1430 915
\pinlabel $(f)$ at 30 425
\pinlabel $(g)$ at 389 425
\pinlabel $(h)$ at 705 425
\pinlabel $(i)$ at 1055 425
\pinlabel $(j)$ at 1415 425

\pinlabel $h$ at 170 555
\pinlabel $h$ at 479 555
\pinlabel $h$ at 822 555
\pinlabel $h$ at 1205 555
\pinlabel $h$ at 1605 555

\pinlabel $h$ at 209 65
\pinlabel $h$ at 553 65
\pinlabel $h$ at 850 65
\pinlabel $h$ at 1180 65
\pinlabel $h$ at 1525 65

\endlabellist 
\begin{center}
\includegraphics[height = 8cm]{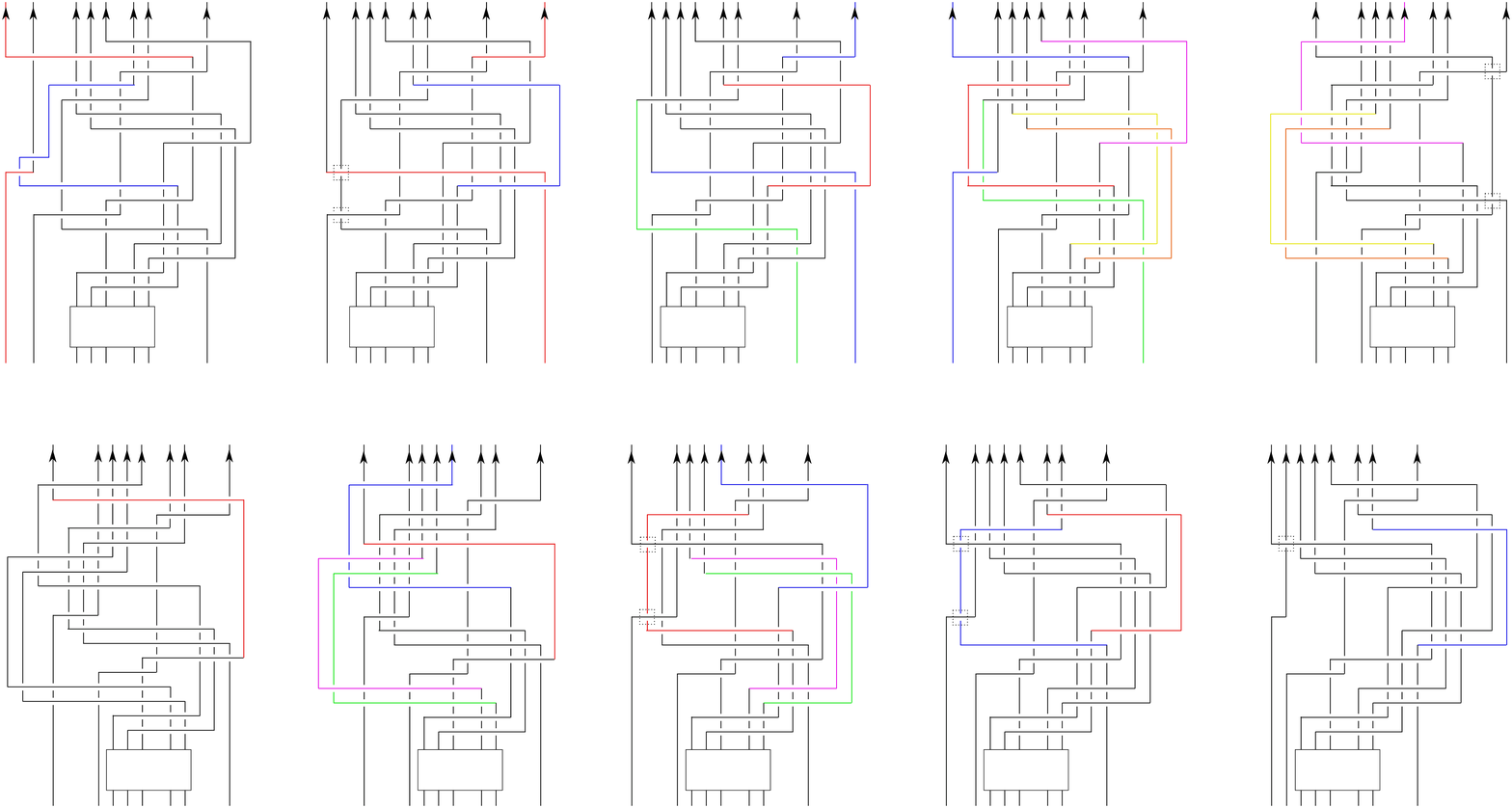}
\caption{\quad Figure \ref{fig:grids2}.a shows a rectangular braid diagram for $T_{hw_1}$. After a sequence of braid isotopies and exchange moves, we obtain the braid in \ref{fig:grids2}.j, which can be negatively destabilized. }
\label{fig:grids2}
\end{center}
\end{figure}

According to Proposition \ref{prop:EH}, $T_{w_1}$ is transversely isotopic to a braid which can be negatively destabilized. Figure \ref{fig:grids2} shows a sequence of transverse braid moves which demonstrates that the same is true of $T_{hw_1}$ for any word $h\in B_8$ in the generators $\sigma_3,\dots, \sigma_6$. The braid in Figure \ref{fig:grids2}.b is obtained from that in \ref{fig:grids2}.a by isotoping the red and blue arcs as shown. The braid in \ref{fig:grids2}.c is related to that in \ref{fig:grids2}.b by an exchange move at the circled crossings in \ref{fig:grids2}.b. The braid in \ref{fig:grids2}.d is obtained from that in \ref{fig:grids2}.c by isotopy of the red, blue and green arcs. The braid in \ref{fig:grids2}.e is obtained from that in \ref{fig:grids2}.d after the indicated isotopy of the yellow, orange and purple arcs. An exchange move at the circled crossings in \ref{fig:grids2}.e produces the braid in \ref{fig:grids2}.f. The braid in \ref{fig:grids2}.g is obtained from that in \ref{fig:grids2}.f by isotoping the red arc as shown. The braid in \ref{fig:grids2}.h is obtained from that in \ref{fig:grids2}.g after an isotopy of the blue, green and purple arcs as shown. An exchange move at the circled crossings in \ref{fig:grids2}.h, followed by the indicated isotopy of the red arc produces the braid in \ref{fig:grids2}.i. Finally, the braid in \ref{fig:grids2}.j is obtained from that in \ref{fig:grids2}.i by an exchange move at the circled crossings in \ref{fig:grids2}.i, followed by the indicated isotopy of the blue arc. Note that the braid in \ref{fig:grids2}.j may be negatively destabilized at the circled crossing. The essential point here is that the region of the braid in \ref{fig:grids2}.a corresponding to the word $h$ is not affected by this combination of isotopies and exchange moves.

To sum up: since $\widehat\theta(T_{w_2}) \neq 0$ (see \cite{not}), we have proven that for any $h\in B_8$ which is 1) a word in the generators $\sigma_3,\dots,\sigma_6$ and for which 2) $\widehat\theta(T_{h}) \neq 0$, it is the case that $\widehat\theta(T_{hw_1})=0$ while $\widehat\theta(T_{hw_2})\neq 0$. It follows that the transverse braids $T_{hw_1}$ and $T_{hw_2}$ are not transversely isotopic though they are topologically isotopic. If, in addition, 3) $h$ is such that the two strands of $T_{hw_1}$ which cross according to the string $\sigma_1^2$ belong to the same component of $T_{hw_1}$, then $\mathcal{SL}(T_{hw_1}) = \mathcal{SL}(T_{hw_2})$; that is, the topological link type represented by $T_{hw_1}$ is transversely non-simple. 


There are infinitely many choices of $h$ which meet criteria 1) - 3) above. In order to give such an $h$, we first prove the following. 

\begin{lemma}
\label{lem:translation}
For $1\leq j\leq k$ and $0\leq l\leq k-j$, consider the map $\psi_{j,k,l}:B_j \rightarrow B_k$ which sends $\sigma_i$ to $\sigma_{i+l}$. If $g$ is a word in $B_j$ for which $\widehat\theta(T_{g})\neq 0$, and $h = \psi_{j,k,l}(g)$, then $\widehat\theta(T_{h}) \neq 0$ as well. 
\end{lemma}

\begin{figure}[!htbp]
\labellist 
\hair 2pt 
\small
\pinlabel $g$ at 53 100
\pinlabel $g$ at 289 100
\pinlabel $T_g$ at 52 0
\pinlabel $T_{\psi_{j,k,l}(g)}$ at 288 0

\tiny
\pinlabel $j$ at 52 215
\pinlabel $j$ at 288 215
\pinlabel $l$ at 210 215
\pinlabel $k-j-l$ at 370 215

\endlabellist 
\begin{center}
\includegraphics[height = 3.1cm]{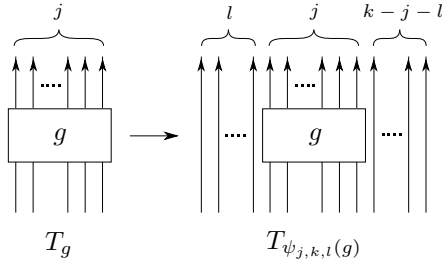}
\caption{\quad On the left, the $j$-braid $T_g$. On the right, the $k$-braid $T_{\psi_{j,k,l}(g)}$. }
\label{fig:map}
\end{center}
\end{figure}

See Figure \ref{fig:map} for a pictorial depiction of the map $\psi_{j,k,l}$.

\begin{proof}[Proof of Lemma \ref{lem:translation}]
If $h=\psi_{j,k,l}(g)$, then the braid $T_h$ is easily seen to be connected sum of $T_g$ with the trivial braids $I_{k-j-l+1}$ and $I_{l+1}$. We know, from the proof of Corollary \ref{cor:qp}, that $\widehat\theta(I_n) \neq 0$ for any $n\geq 1$. Lemma \ref{lem:translation} therefore follows from Theorem \ref{thm:connectedsum}.
\end{proof}

It follows from Corollary \ref{cor:qp} and from Lemma \ref{lem:translation} that $h=\psi_{4,8,3}(g)$ satisfies criteria 1) and 2) above as long as $T_g$ is a quasipositive 4-braid. Let $h=\psi_{4,8,3}(g)$ for $$g=\sigma_3\sigma_2\sigma_3\sigma_1\sigma_2\sigma_3.$$ It is easy to check that $h^n$ also satisfies criterion 3) (as well as criteria 1) and 2), of course) for all $n\geq 0$.
 
\begin{corollary}
\label{cor:ex}
The topological link types represented by $T_{h^nw_1}$ are transversely non-simple for all $n\geq 0$. When $n$ is even, $T_{h^nw_1}$ is a knot; otherwise, $T_{h^nw_1}$ is a 3-component link.
\end{corollary}

\begin{figure}[!htbp]
\labellist 
\hair 2pt 
\small
\pinlabel $T_{w_1}$ at 140 20
\pinlabel $T_{h^{2m}w_1}$ at 468 20
\pinlabel $w_1$ at 140 262
\pinlabel $w_1$ at 468 262
\pinlabel $w_1$ at 796 262
\pinlabel $w_1$ at 1116 262
\pinlabel $h^{2m}$ at 468 150
\pinlabel $m$ at 826 146
\pinlabel $-\frac{1}{m}$ at 1296 55

\tiny
\pinlabel $=$ at 630 190
\pinlabel $=$ at 950 190

\endlabellist 
\begin{center}
\includegraphics[height = 3.2cm]{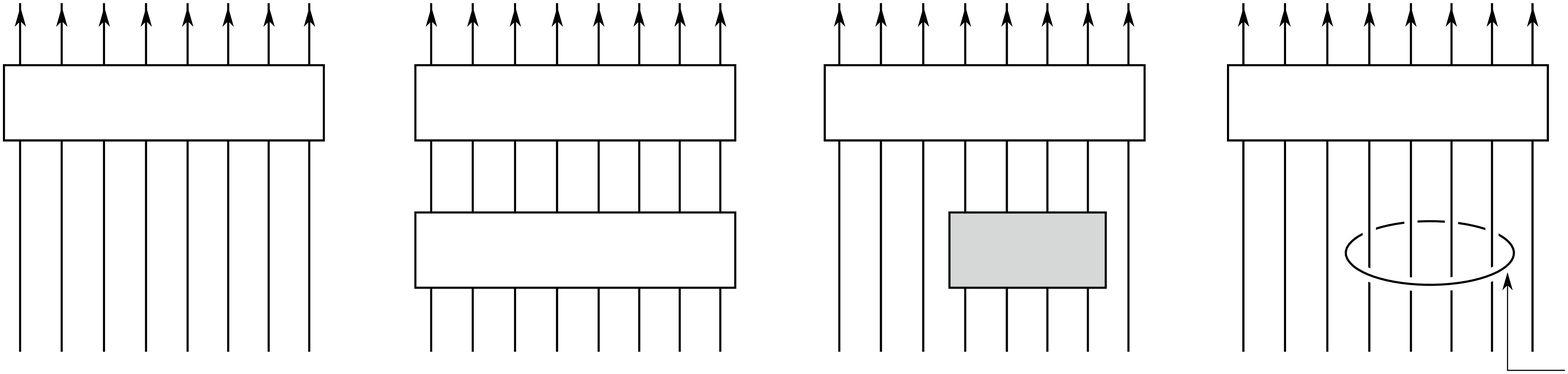}
\caption{\quad $T_{h^{2m}w_1}$ is obtained from $T_{w_1}$ by performing $m$ positive full twists of strands 4 - 7, or, equivalently, by performing $-1/m$ surgery on the unknot shown on the right.}
\label{fig:braid2}
\end{center}
\end{figure}

Below, we prove that most of the links in Corollary \ref{cor:ex} are prime. Note that $T_{h^nw_1}$ is obtained from $T_{w_1}$ by performing $n$ positive half twists of strands 4 - 7 in the region of $T_{w_1}$ where we would insert the word $h^n$. For $n=2m$, this amounts to adding $m$ positive full twists, which can also be accomplished by performing $-1/m$ surgery on an unknot $U$ encircling strands 4 - 7 of $T_{w_1}$ in the corresponding region. See Figure \ref{fig:braid2}.

A \emph{SnapPea} computation \cite{weeks} combined with the Inverse Function Theorem test described in Moser's thesis \cite{moserh} shows that the complement of the link $T_{w_1} \cup U$ is hyperbolic. To be specific, \emph{SnapPea} finds a triangulation of this link complement by ideal tetrahedra and computes an approximate solution to the gluing equations. Moser's test then confirms, using this approximate solution, that an exact solution exists. 

Thurston's celebrated Dehn Surgery Theorem then implies that all but finitely many Dehn fillings of the boundary component corresponding to $U$ are hyperbolic as well \cite{th1}. In turn, this implies that the link $T_{h^{2m}w_1}$ is hyperbolic, and, hence, prime for all but finitely many $m$. This argument can be repeated to show that the links $T_{h^{2m+1}w_1}$ are also prime for all but finitely many $m$. The lemma below sums this up.

\begin{lemma}
The links $T_{h^nw_1}$ are prime for all but finitely many values of $n$.
\end{lemma}

\bibliographystyle{hplain.bst}
\bibliography{References}

\begin{thebibliography}{10}

\bibitem{bev}
K.~Baker, J.B. Etnyre, and J.~Van Horn-Morris.
\newblock Fibered transverse knots and the bennequin bound.
\newblock 2008, math.GT/0803.0758.

\bibitem{bald3}
J.~A. Baldwin.
\newblock Comultiplicativity of the {O}zsv{\'a}th-{S}zab{\'o} contact
  invariant.
\newblock {\em Math. Res. Lett.}, 15(2):273--287, 2008.

\bibitem{benn}
D.~Bennequin.
\newblock Entrelacements et {\'e}quations de {P}faff.
\newblock {\em Ast{\'e}risque}, 107-108:87--161, 1983.

\bibitem{bm2}
J.S. Birman and W.W. Menasco.
\newblock Studying links via closed braids {IV}: composite links and split
  links.
\newblock {\em Inv. Math.}, 102(1):115--139, 1990.

\bibitem{bm3}
J.S. Birman and W.W. Menasco.
\newblock {Stabilization in the braid group II: Transversal simplicity of
  knots}.
\newblock {\em Geom. Topol.}, 10:1425--1452, 2006.

\bibitem{bw}
J.S. Birman and R.F. Williams.
\newblock {Knotted periodic orbits in dynamical systems - I: Lorenz's
  equations}.
\newblock {\em Topology}, 22(1):47--82, 1983.

\bibitem{yasha5}
Y.~Eliashberg.
\newblock Legendrian and transversal knots in tight contact 3-manifolds.
\newblock In {\em Topological methods in modern mathematics}, pages 171--193.
  {Publish or Perish}, 1993.

\bibitem{efm}
J.~Epstein, D.~Fuchs, and M.~Meyer.
\newblock {Chekanov-Eliashberg invariants and transverse approximations of
  Legendrian knots}.
\newblock {\em Pac. J. Math.}, 201(1):89--106, 2001.

\bibitem{et4}
J.~B. Etnyre.
\newblock Transversal torus knots.
\newblock {\em Geom. Topol.}, 3:253--268, 1999.

\bibitem{EH2}
J.~B. Etnyre and K.~Honda.
\newblock Knots and contact geometry {I}: torus knots and the figure eight
  knot.
\newblock {\em J. Symp. Geom.}, 1(1):63--120, 2001.

\bibitem{EH5}
J.~B. Etnyre and K.~Honda.
\newblock On connected sums and legendrian knots.
\newblock {\em Adv. Math.}, 179(1):59--74, 2003.

\bibitem{EH4}
J.~B. Etnyre and K.~Honda.
\newblock Cabling and transverse simplicity.
\newblock {\em Ann. of Math.}, 162(3):1305--1333, 2005.

\bibitem{keiko}
K.~Kawamuro.
\newblock Connect sum and transversly non-simple knots.
\newblock {\em Math. Proc. Cambridge Philos. Soc., to appear}, 2008.

\bibitem{kng}
T.~Khandhawit and L.~Ng.
\newblock A family of transversely nonsimple knots.
\newblock {\em Algebr. Geom. Topol.}, 10(1):293--314, 2010.

\bibitem{lossz}
P.~Lisca, P.~Ozsv{\'a}th, A.~Stipsicz, and Z.~Szab{\'o}.
\newblock Heegaard {F}loer invariants of {L}egendrian knots in contact
  three-manifolds.
\newblock 2008, math.SG/0802.0628.

\bibitem{mos}
C.~Manolescu, P.~Ozsv{\'a}th, and S.~Sarkar.
\newblock A combinatorial description of knot {F}loer homology.
\newblock {\em Annals of Mathematics}, 169:633--660, 2009.

\bibitem{most}
C.~Manolescu, P.~Ozsv{\'a}th, Z.~Szab{\'o}, and D.~Thurston.
\newblock On combinatorial link {F}loer homology.
\newblock {\em Geom. Topol.}, 11:2339--2412, 2007.

\bibitem{mm}
W.~W. Menasco and H.~Matsuda.
\newblock An addendum on iterated torus knots (appendix).
\newblock 2006, math.GT/0610566.

\bibitem{moserh}
H.~H. Moser.
\newblock {\em Proving a manifold to be hyperbolic once it has been
  approximated to be so}.
\newblock PhD thesis, Columbia University, 2005.

\bibitem{not}
L.~Ng, P.~Ozsv{\'a}th, and D.~Thurston.
\newblock Transverse knots distinguished by knot {F}loer homology.
\newblock {\em J. Symp. Geom.}, 6(4):461--490, 2008.

\bibitem{ngt}
L.~Ng and D.~Thurston.
\newblock Grid diagrams, braids, and contact geometry.
\newblock pages 120--136, 2009.

\bibitem{osh}
S.~Orevkov and V.~Shevchishin.
\newblock Markov {T}heorem for {T}ransverse {L}inks.
\newblock {\em J. Knot Theory Ram.}, 12(7):905--913, 2003.

\bibitem{ost}
P.~Ozsv{\'a}th and A.~Stipsicz.
\newblock Contact surgeries and the transverse invariant in knot floer
  homology.
\newblock 2008, math.GT/0803.1252.

\bibitem{osz19}
P.~Ozsv{\'a}th and Z.~Szab{\'o}.
\newblock Holomorphic disks, link invariants, and the multi-variable
  {A}lexander polynomial.
\newblock {\em Algebr. Geom. Topol.}, 8:615--692, 2008.

\bibitem{oszt}
P.~Ozsv{\'a}th, Z.~Szab{\'o}, and D.~Thurston.
\newblock Legendrian knots, transverse knots, and combinatorial {F}loer
  homology.
\newblock {\em Geom. Topol.}, 12:941--980, 2008.

\bibitem{th1}
W.~Thurston.
\newblock {\em The geometry and topology of three-manifolds}.
\newblock Princeton, 1979.

\bibitem{vv}
D.~S. Vela-Vick.
\newblock On the transverse invariant for bindings of open books.
\newblock 2009, math.SG/0806.1729.

\bibitem{vera}
V.~V{\'e}rtesi.
\newblock Transversely non-simple knots.
\newblock {\em Algebr. Geom. Topol.}, 8:1481--1498, 2008.

\bibitem{weeks}
J.~Weeks.
\newblock Snap{P}ea.
\newblock { \tt http://www.geometrygames.org/Snap{P}ea/index.html}.

\bibitem{wrinkle}
N.~Wrinkle.
\newblock The {M}arkov theorem for transverse knots.
\newblock 2002, math.GT/0202055.

\end{thebibliography}

\end{document}